\documentclass[reqno,centertags]{amsart}
\usepackage{amsmath}
\usepackage{amssymb}
\usepackage{amsfonts}

\newtheorem{theorem}{Theorem}
\theoremstyle{plain}

\newtheorem{problem}{Problem}

\numberwithin{equation}{section}

\input{tcilatex}

\begin{document}
\title[]{Anchor rings of finite type Gauss map in the Euclidean 3-space}
\author{Hassan Al-Zoubi}
\address{Department of Mathematics, Al-Zaytoonah University of Jordan, P.O.
Box 130, Amman, Jordan 11733}
\email{dr.hassanz@zuj.edu.jo}
\author{Nancy Al-Ramahi}
\address{Department of Basic Sciences, Al-Zaytoonah University of Jordan, P.O.
Box 130, Amman, Jordan 11733}
\email{N.Alramahi@zuj.edu.jo}
\date{}

\begin{abstract}
In this article, we continue the classification of finite type Gauss map surfaces in the Euclidean 3-space $\mathbb{E}^{3}$ with respect to the first fundamental form by studying a subclass of tubes, namely the anchor rings. We show that anchor rings are of infinite type Gauss map.
\end{abstract}

\maketitle

\section{Introduction}
The theory of Gauss map of a surface in a Euclidean space has been investigated from the various viewpoints by many differential geometers \cite{C10, C11, C12, C13, D5, K1, K2, K3, K4, K5, K6, K7, R2}, and it has been a topic of active research. In this paper, we will be concerned with the theory of finite type Gauss map of surfaces in the Euclidean 3-space. Related to this subject many interesting results have been found in \cite{B2, B3, B5, B6, D4}.

In 1983 B.-Y. Chen introduced the notion of Euclidean immersions of finite type \cite{C3}. A surface $S$ is said to be of finite type corresponding to the first fundamental form $I$, or briefly of finite $I$-type, if the position vector $\boldsymbol{x}$ of $S$ can be written as a finite sum of nonconstant eigenvectors of the Laplacian $\Delta ^{I}$, that is,
\begin{equation}  \label{1}
\boldsymbol{x}=\boldsymbol{c}+\sum_{i=1}^{k}\boldsymbol{x}_{i},\quad
\Delta ^{I}\boldsymbol{x}_{i}=\lambda _{i}\,\boldsymbol{x}_{i},\quad
i=1,\dotsc ,k,
\end{equation}%
where $\boldsymbol{c}$ is a fixed vector and $\lambda _{1},\lambda_{2},\dotsc ,\lambda _{k}$ are eigenvalues of the operator $\Delta^{I}$. In particular, if all eigenvalues $\lambda _{1},\lambda _{2},\dotsc ,\lambda_{k}$ are mutually distinct, then $S$ is said to be of finite $I$-type $k$. When $\lambda _{i}=0$ for some $i=1,\dotsc ,k$, then $S$ is said to be of finite null $I$-type $k$. Otherwise, $S$ is said to be of infinite type.

In general when $S$ is of finite type $k$, it follows from (\ref{1}) that there exists a monic polynomial, say $F(x)\neq 0,$ such that $F(\Delta ^{I})(\boldsymbol{x}-\boldsymbol{c})=\mathbf{0}.$ Suppose that $F(x)=x^{k}+\sigma_{1}x^{k-1}+...+\sigma _{k-1}x+\sigma _{k},$ then coefficients $\sigma _{i}$ are given by

\begin{eqnarray}
\sigma _{1} &=&-(\lambda _{1}+\lambda _{2}+...+\lambda _{k}),  \notag \\
\sigma _{2} &=&(\lambda _{1}\lambda _{2}+\lambda _{1}\lambda_{3}+...+\lambda_{1}\lambda _{k}+\lambda _{2}\lambda _{3}+...+\lambda _{2}\lambda
_{k}+...+\lambda _{k-1}\lambda _{k}),  \notag \\
\sigma _{3} &=&-(\lambda _{1}\lambda _{2}\lambda _{3}+...+\lambda_{k-2}\lambda _{k-1}\lambda _{k}),  \notag \\
&&.............................................  \notag \\
\sigma _{k} &=&(-1)^{k}\lambda _{1}\lambda _{2}...\lambda _{k}.  \notag
\end{eqnarray}
Therefore, the position vector $\boldsymbol{x}$ satisfies the following equation (see \cite{C7})

\begin{equation*}  \label{2}
(\Delta^{I})^{k}\boldsymbol{x}+\sigma_{1}(\Delta^{I})^{k-1}\boldsymbol{x}+...+\sigma_{k}( \boldsymbol{x}-\boldsymbol{x}_{0})=\boldsymbol{0}.
\end{equation*}

Very little is known about surfaces of finite type in the Euclidean 3-space $\mathbb{E}^{3}$. Actually, so far, the only known surfaces of finite $I$-type in $\mathbb{E}^{3}$ are the minimal surfaces, the circular cylinders, and the spheres. So in \cite{C2} B.-Y. Chen mentions the following problem
\begin{problem}  \label{(p1)}
classify all surfaces of finite Chen $I$-type in $\mathbb{E}^{3}$.
\end{problem}

Many authors have been interested and studied this problem by investigating special classes of surfaces. More precisely, starting from the late 1980s, the above problem was solved for the class of spiral surfaces, tubes, ruled surfaces, quadric surfaces and the compact and noncompact cyclides of Dubin, see \cite{B1, C4, C5, C6, D1, D2}. Meanwhile, this problem still not solved yet for the class of surfaces of revolution, translation surfaces, as well as helicoidal surfaces.

We consider the surface $S$ in $\mathbb{E}^{3}$. The map $\boldsymbol{n} : S \rightarrow M^{2}$ which sends each point of $S$ to the unit normal vector to $S$ at the point is called the Gauss map of the surface $S$, where $M^{2}$ is the unit sphere in $\mathbb{E}^{3}$ centered at the origin. In this regard, the above problem can be generalized by studying surfaces whose Gauss map $\boldsymbol{n}$ is of finite type with respect to the first fundamental form (see \cite{C9}), so specifically in the Euclidean 3-space, we pose the following problem
\begin{problem}  \label{p2}
Classify all surfaces with finite type Gauss map in $\mathbb{E}^{3}$.
\end{problem}

Regarding to this, problem \ref{p2} was solved for the class of spiral surfaces \cite{B1}, cyclides of Dupin \cite{B3}, ruled surfaces and tubes \cite{B4}. However, surfaces of revolution, translation surfaces, cones, quadric surfaces, as well as helicoidal surfaces, the classification of its finite type Gauss map is not known yet.

\section{Preliminaries}

Let $\boldsymbol{x} = \boldsymbol{x}(u^{1}, u^{2})$ be a regular parametric representation of $S$ in $\mathbb{E}^{3}$. For a sufficient differentiable function $f(u^{1}, u^{2})$ the second Beltrami-Laplace operator with respect to the first fundamental form $I = g_{ij}du^{i}du^{j}$ of $S$ is defined by

\begin{equation}  \label{3}
\Delta ^{I}f =-\frac{1}{\sqrt{g}}\Big(\sqrt{g}g^{ij}f_{i}\Big)_{j},
\end{equation}
where $g:=det(g_{ij})$ and $g^{ij}$ denote the components of the inverse tensor of $g_{ij}$.
Applying (\ref{3}) for the position vector $\boldsymbol{x}$, we have the following well-known formula
\begin{equation}  \label{4}
\Delta ^{I}\boldsymbol{x}=-2H\boldsymbol{n},
\end{equation}
where $H$ denotes the mean curvature of $S$.
From (\ref{4}) we know the following two facts \cite{T1}

\begin{itemize}
\item $S$ is minimal if and only if all coordinate functions of $\boldsymbol{x}$ are eigenfunctions of $\Delta ^{I}$ with eigenvalue 0.

\item $S$ lies in an ordinary sphere $M^{2}$ if and only if all coordinate functions of $\boldsymbol{x}$ are eigenfunctions of $\Delta ^{I}$ with a fixed nonzero eigenvalue.
\end{itemize}

By applying (\ref{3}) for the normal vector $\boldsymbol{n}$ we get \cite{S1}
\begin{equation*}  \label{5}
\Delta ^{I}\boldsymbol{n}=grad^{I}2H +(4H^{2}-2K)\boldsymbol{n},
\end{equation*}
where $K$ denotes the Gaussian curvature of $S$. Up to now, the only known surfaces of finite type Gauss map are the spheres, the minimal surfaces, and the circular cylinders.
In the present paper, we mainly focus on problem \ref{p2} by studying a subclass of tubes in $\mathbb{E}^{3}$, namely anchor rings.

\section{Anchor rings in $\mathbb{E}^{3}$}

First we define tubes in the Euclidean 3-space. Let $C: \boldsymbol{\alpha} =\boldsymbol{\alpha}(\emph{t})$, $\mathit{t}\epsilon (a,b)$ be a regular unit speed curve of finite length which is topologically imbedded in $\mathbb{E}^{3}$. The total space $N_{\boldsymbol{\alpha}}$ of the normal bundle of $\boldsymbol{\alpha}((a, b))$ in $\mathbb{E}^{3}$ is naturally diffeomorphic to the direct product $(a,b)\times \mathbb{E}^{2}$ via the translation along $\boldsymbol{\alpha}$ with respect to the induced normal connection. For a sufficiently small $r>0$ the tube of radius $r$ about the
curve $\boldsymbol{\alpha}$ is the set:
\begin{equation}
T_{r}( \boldsymbol{\alpha})=\{exp_{\boldsymbol{\alpha}(t)}\boldsymbol{u}\mid \boldsymbol{u}\in N_{\boldsymbol{\alpha}} , \ \ \parallel \boldsymbol{u}%
\parallel =r,\ \ t\in(a,b)\}.  \notag
\end{equation}
Assume that ${\mathbf{t}, \mathbf{h}, \mathbf{b}}$ is the Frenet frame and that $\kappa$ is the curvature of the unit speed curve $\boldsymbol{\alpha} =
\boldsymbol{\alpha}(\emph{t})$. For a small real number $r$ satisfies $0 < r< min\frac{1}{|\kappa|}$, the tube $T_{r}( \boldsymbol{\alpha})$ is a smooth surface in $\mathbb{E}^{3}$ \cite{R1}. Then, a parametric representation of the tube $T_{r}( \boldsymbol{\alpha})$ is given by

\begin{equation*}  \label{eq6}
\digamma:\boldsymbol{x}(t,\varphi)= \boldsymbol{\alpha}+ r \cos\varphi\mathbf{h}+ r \sin\varphi\mathbf{b}.
\end{equation*}

It is easily verified that the first fundamental form of $\digamma$ is given by
\begin{align*}
I &= \big(\delta^{2}+r^{2}\tau^{2}\big)dt^{2} + 2r^{2}\tau dtd\varphi+r^{2}d\varphi^{2},
\end{align*}
where $\delta: = (1-r\kappa\cos\varphi)$ and $\tau$ is the torsian of the curve $\boldsymbol{\alpha}$.
The Beltrami operator corresponding to the first fundamental form of $\digamma$ can be expressed as follows

\begin{eqnarray*}  \label{eq8}
\Delta ^{I} &=&-\frac{1}{(\delta)^{3}}\Bigg[\delta\frac{\partial^{2}}{\partial t^{2}}-2\tau \delta\frac{\partial ^{2}}{\partial t\partial \varphi}%
+\frac{\delta}{r^{2}}(r^{2}\tau ^{2}+\delta^{2})\frac{\partial^{2}}{\partial \varphi^{2}}\\
&&+r\beta\frac{\partial}{\partial t}-\frac{\kappa \delta^{2}\sin\varphi}{r}\frac{\partial }{\partial \varphi }\Bigg] , \notag
\end{eqnarray*}
where $\beta: =\kappa \acute{} \cos \varphi +\kappa \tau \sin \varphi $ and $\acute{}:=\frac{d}{dt}$.

Now, we define an anchor ring in the Euclidean 3-space. A tube in $\mathbb{E}^{3}$ is called an anchor ring if the curve $C$ is a plane circle (or is an open portion of a plane circle). In this case, the torsian $\tau$ of $\alpha$ vanishes identically and the curvature $\kappa$ of $\alpha$ is a nonzero constant. Then, the position vector $\boldsymbol{x}$ of the anchor ring can be expressed as  \cite{A3,A5}

\begin{equation}  \label{eq9}
\digamma:\boldsymbol{x}(t,\varphi)= \{\gamma \cos\varphi,\gamma \sin\varphi, r\sin t \},
\end{equation}

\begin{center}
$a > r,\ \ a \epsilon \mathbb{R},$
\end{center}
where $\gamma :=a+r\cos t$. Then the first fundamental form of (\ref{eq9}) is

\begin{equation*}
I=r^{2}dt^{2}+\gamma^{2}d\varphi ^{2}.
\end{equation*}
Hence, the Beltrami operator is given by

\begin{equation}  \label{eq10}
\Delta^{I}= -\frac{1}{\gamma^{^{2}}}\frac{\partial^{2}}{\partial \varphi^{2}}+\frac{\sin t}{r\gamma}\frac{\partial}{\partial t} - \frac{1}{r^{^{2}}}\frac{\partial^{2}}{\partial t^{2}}.
\end{equation}

Denoting by $\boldsymbol{n}$ the Gauss map of $\digamma$, then we have
\begin{equation*}
\boldsymbol{n}=\{-\cos t \cos \varphi, -\cos t \sin\varphi, -\sin t\}.
\end{equation*}

Let $n_{3}$ be the third coordinate function of $\boldsymbol{n}$. By virtue of (\ref{eq10}), one can find

\begin{equation}  \label{eq13}
\Delta^{I}n_{3} = -\frac{\sin t}{r}\Bigg[\frac{\cos t}{\gamma} +\frac{1}{r}\Bigg].
\end{equation}

Moreover, by direct computation, we obtain

\begin{equation*}
(\Delta^{I})^{2}n_{3} =-\frac{\sin t}{r^{4}}-\frac{5\sin t \cos t}{r^{3}\gamma}-\frac{\sin^{3}t}{r^{2}\gamma^{2}}+\frac{2\cos^{2}t \sin t}{r^{2}\gamma^{2}}
-\frac{3\sin^{3}t \cos t}{r\gamma^{3}},
\end{equation*}
which can be rewritten as
\begin{equation}  \label{eq14}
(\Delta^{I})^{2}n_{3} =-\frac{3\sin^{3}t \cos t}{r\gamma^{3}}+\frac{1}{\gamma^{2}}F_{2}(\sin t,\cos t),
\end{equation}
where $F_{2}(\sin t,\cos t)$ is a polynomial in $\sin t$ and $\cos t$ of degree 3. Applying relation (\ref{eq10}) on (\ref{eq14}) gives

\begin{equation}  \label{eq15}
(\Delta^{I})^{3}n_{3} =-\frac{45\sin^{5}t \cos t}{r\gamma^{5}}+\frac{1}{\gamma^{4}}F_{3}(\sin t,\cos t),
\end{equation}
where $F_{3}(\sin t,\cos t)$ is a polynomial in $\sin t$ and $\cos t$ of degree 5. For each integer $k>0$ it can be easily seen that

\begin{equation}  \label{eq16}
\Delta^{I}\frac{\sin^{k}t \cos t}{r\gamma^{k}} = \lambda_{k}\frac{\sin^{2k-1}t \cos t}{r\gamma^{2k-1}}+\frac{1}{\gamma^{2k-2}}Q_{k}(\sin t,\cos t),
\end{equation}
where $Q_{k}(\sin t,\cos t)$ is a polynomial in $\sin t$ and $\cos t$ of degree $2k-1$ and $\lambda_{k}=\overset{k}{\underset{j=1}{\prod }}(2j-1)(2j-3)$.
Therefore, one can find

\begin{equation}  \label{eq17}
(\Delta^{I})^{m}n_{3} = \lambda_{m}\frac{\sin^{2m-1}t \cos t}{r\gamma^{2m-1}}+\frac{1}{\gamma^{2m-2}}F_{m}(\sin t,\cos t),
\end{equation}%
where $F_{m}(\sin t,\cos t)$ is a polynomial in $\sin t$ and $\cos t$ of degree $2m-1$.

Notice that $\lambda_{k}\neq 0$, for each natural number $k$. Now, if the Gauss map $\boldsymbol{n}$ is of finite type, then there exist real numbers, $c_{1}, c_{2}, ... , c_{m}$ such that

\begin{equation}  \label{eq18}
(\Delta ^{I})^{m}\boldsymbol{n}+c_{1}(\Delta ^{I})^{m-1}\boldsymbol{n}+...+c_{m-1}\Delta ^{I}\boldsymbol{n}+c_{m}\boldsymbol{n}=\mathbf{0}.
\end{equation}

Since $n_{3}=-\sin t $ is the third coordinate of $\boldsymbol{n}$, one gets

\begin{equation}  \label{eq17}
(\Delta ^{I})^{m}n_{3}+c_{1}(\Delta ^{I})^{m-1}n_{3}+...+c_{m-1}\Delta^{I}n_{3}+c_{m}n_{3}= 0.
\end{equation}

From (\ref{eq13}-\ref{eq15}), (\ref{eq17}) and (\ref{eq18}) we obtain that

\begin{eqnarray*}
\lambda_{m}\frac{\sin^{2m-1}t \cos t}{r\gamma^{2m-1}}+\frac{1}{\gamma^{2m-2}}F_{m}(\sin t,\cos t)&&  \notag \\
+c_{1}\lambda_{m-1}\frac{\sin^{2m-3}t \cos t}{r\gamma^{2m-3}}+c_{1}\frac{1}{\gamma^{2m-4}}F_{m-1}(\sin t,\cos t)&&  \notag \\
+...+c_{m-1}\frac{\sin t}{r}\Bigg[\frac{\cos t}{\gamma} +\frac{1}{r}\Bigg]+c_{m}\sin t &=& 0
\end{eqnarray*}
which can be rewritten as

\begin{equation}  \label{eq18}
\lambda_{m}\frac{\sin^{2m-1}t \cos t}{r\gamma}+R(\cos t, \sin t)=0,
\end{equation}%
where $R(\cos t, \sin t)$ is a polynomial in $\cos t$ and $\sin t$ of degree $2m-1$.

This is impossible for any $m \geq 1$ since $\lambda_{m} \neq 0$. Consequently, we have the following

\begin{theorem}
\label{C1.1} Every anchor ring in the Euclidean 3-space is of infinite type.
\end{theorem}




\end{document}